\documentclass[11pt]{article}
\usepackage{latexsym}
\usepackage{amsfonts}
\usepackage{mathrsfs,amsthm,amsmath}
\usepackage{todonotes}
\usepackage{color}
\newtheorem{theorem}{Theorem}[section]
\newtheorem{proposition}[theorem]{Proposition}
\newtheorem{definition}{Definition}[section]

\newtheorem{lemma}{Lemma}[section]
\newtheorem{remark}{Remark}[section]
\newtheorem{example}{Example}[section]
\newtheorem{condition}{Condition}[section]
\begin{document}
\title{Large and moderate deviations for Gaussian neural networks}
\author{Claudio Macci\thanks{Address: Dipartimento di Matematica, Università di Roma Tor Vergata,
		Via della Ricerca Scientifica, I-00133 Rome, Italy. e-mail: \texttt{macci@mat.uniroma2.it}}
	\and Barbara Pacchiarotti\thanks{Address: Dipartimento di Matematica, Università di Roma Tor Vergata,
		Via della Ricerca Scientifica, I-00133 Rome, Italy. e-mail: \texttt{pacchiar@mat.uniroma2.it}}
	\and Giovanni Luca Torrisi\thanks{Address:
        Istituto per le Applicazioni del Calcolo, Consiglio Nazionale delle Ricerche, Via dei Taurini 19, 00185 
        Rome, Italy. email: \texttt{giovanniluca.torrisi@cnr.it}}}
\maketitle

\begin{abstract}
We prove large and moderate deviations for the output of Gaussian fully connected neural networks. The main achievements concern deep neural networks (i.e.,  when the model has more than one hidden layer) and hold for bounded and continuous pre-activation functions. However,  for deep neural networks fed by a single input,  we have results even if the pre-activation is ReLU.  When the network is shallow (i.e., there is exactly one hidden layer) the large and moderate principles hold for quite general pre-activation functions.
\end{abstract}

\section{Introduction}

In the last decade neural networks have been successfully exploited to solve a variety of practical problems,  ranging from computer vision and speech recognition to feature extraction \cite{GPMX,  LBH}.  
This has stimulated new mathematical research in different fields such as probability and statistics, with the final goal to better understand how neural networks work and how to make them more efficient \cite{ADFST,  BGRS, BKLW,  BT,  BFF23, CMSV,  EMS21, FFP, FHMNP, Hanin, 
Hirsch2024, JungLeeLeeYang, SS, Vogel, zavatoneveth2021}.  Indeed,  despite their profound engineering
success,  a comprehensive understanding of the intrinsic working mechanism of neural networks is still lacking.  In particular,  the analysis of {\em deep} neural networks is very challenging due
to the recursive and nonlinear structure of the models.

Neural networks are parametrized families of functions which are typically used in statistical learning to estimate an unknown function $f$.  In practice,  one first fixes the network architecture,  specifying in this way 
the family of parametric functions,  and then looks for an approximation of the target function $f$,  within the specified family, on the basis of a given training set of data
\cite{RYH}.

In this paper,  we focus on the class of {\em fully connected neural networks} which are formally defined as follows.  Fix a positive integer $L\geq 1$, $L+2$ positive integers $n_0, n_1,\ldots,n_{L+1}$ and a function $\sigma:\mathbb R\to\mathbb R$.  A fully connected neural network with depth $L$,  input dimension $n_0$,  output dimension $n_{L+1}$, hidden layer widths $n_1,\ldots,n_L$ and pre-activation function $\sigma:\mathbb R\to\mathbb R$ is a function $$x=(x_{1},\ldots,x_{n_0})\in T\subset\mathbb{R}^{n_0}\mapsto Z^{(L+1)}(x)=(Z_1^{(L+1)}(x),\ldots,Z_{n_{L+1}}^{(L+1)}(x))\in\mathbb{R}^{n_{L+1}}$$ defined by
\begin{equation}\label{eq:Hanin-model}
	\left\{\begin{array}{l}
		Z_i^{(1)}(x)=b_i^{(1)}+\sum_{r=1}^{n_0}W_{ir}^{(1)}x_r\quad(i=1,\ldots,n_1)\\
		Z_i^{(\ell+1)}(x)=b_i^{(\ell+1)}+\sum_{j=1}^{n_\ell}W_{ij}^{(\ell+1)}\sigma(Z_j^{(\ell)}(x))
		  \quad(i=1,\ldots,n_{\ell+1})\ \mbox{for}\ 1\leq\ell\leq L,
	\end{array}\right.
\end{equation}
where $\{b_i^{(\ell)}\}$ (network biases) and $\{W_{ij}^{(\ell)}\}$ (network weights) are the network parameters.
Then,  for fixed $L$ and $n_1,\ldots, n_L$ (i.e., for a fixed network architecture) and given a training dataset\\ $\{(x_\alpha,f(x_\alpha))\}_{\alpha=1,\ldots,m}$,  $m\geq 1$,  one looks for biases and weights such that
$$Z^{(L+1)}(x_\alpha)\approx f(x_\alpha)$$ for inputs $x_\alpha$ not only within,  but also outside 
the training dataset. To this aim,  the usual procedure consists of two steps: $(i)$ One randomly initializes the network parameters $\{b_i^{(\ell)}\}$ and $\{W_{ij}^{(\ell)}\}$; $(ii)$ One optimizes the parameters by minimizing some empirical risk function (such as the squared error).

Therefore, to understand the behavior of a fully connected neural network at the start of training, one studies fully connected neural networks
with random biases and random weights. In most papers in the literature, the network parameters are assumed to be Gaussian distributed (see 
e.g.  \cite{ADFST,  BGRS,BT,  BFF23, FHMNP,  Hanin, RYH}); more specifically, for $C_b\geq 0$ and $C_W>0$, it is supposed that:
\begin{itemize}
	\item the random variables $b_i^{(\ell)}$ are centered Normal distributed with variance $C_b$ ($1\leq \ell\leq L$);
	\item the random variables $W_{ij}^{(\ell+1)}$ are centered Normal distributed with variance $C_W/n_\ell$ ($0\leq\ell\leq L$);
	\item all these random variables are independent.
\end{itemize}
If $L\geq 2$ these neural networks are called Gaussian fully connected {\em deep} neural networks;
if $L=1$ they are called Gaussian fully connected {\em shallow} neural networks.

In literature there is a considerable amount of papers which investigate the asymptotic behavior of fully connected neural networks. Most of them
focus on the distribution approximation of the output in the infinite width limit, i.e.,  when the parameters $n_1,\ldots,n_L$ tend to infinity. 
In the context of shallow neural networks,  a seminal paper in this direction is \cite{Neal}; see  \cite{BGRS,CMSV,EMS21} for related contributions.
An important recent result is 
Theorem 1.2 in \cite{Hanin} which provides, when $T$ is compact, the weak convergence of $((Z_h^{(L+1)}(x))_{x\in T})_{h=1,\ldots,n_{L+1}}$ to
a suitable centered Gaussian field (on $T$),  as $n_1,\ldots,n_L\to\infty$,  for every arbitrarily fixed 
$L\geq 1$ and for continuous and polynomially bounded pre-activation functions.  Quantitative versions of this result,  with respect to different probability metrics, are given in  \cite{ADFST,BT,BFF23,FHMNP}. We emphasize that in \cite{FHMNP} the authors study the Gaussian approximation of the sensitivities of 
the output with respect to the input, i.e., the Gaussian approximation of the mixed directional derivatives of the output with respect to the 
input. Large-width asymptotics of fully connected deep neural networks with biases and weights distributed according to a non-Gaussian stable law 
are investigated in \cite{FFP, JungLeeLeeYang}. The output distribution of a Gaussian fully connected deep neural network, with finite hidden layers widths and a ReLU pre-activation, has been studied in 
\cite{zavatoneveth2021}.

In this paper we are interested in a different kind of asymptotic results,  which are based  on the theory of large deviations (see \cite{DemboZeitouni}).  Such a theory allows to quantify the atypical behavior of the network, and provides probability 
estimates of rare events on an exponential scale. Among the references on large deviations on this topic we recall
\cite{LiSaad} (which concerns a deep neural network model different from \eqref{eq:Hanin-model}), and \cite{Hirsch2024} (whose context is 
different from ours, indeed it concerns the stochastic gradient descent for trained shallow neural networks with quadratic loss). 

In this paper we set $n_\ell=n_\ell(n)$ $(\ell=1,\ldots,L)$,  assume that the sequences $n_\ell$ diverge to infinity 
(as $n\to\infty$)
and, for a suitable normalizing sequence $v_n^*\to\infty$ (see Condition \ref{cond:on-n1-nL}), a
finite set $A$ and $T\equiv\{x_\alpha\}_{\alpha\in A}$ (therefore $T$ is also finite), we 
provide the following main results
for suitable sequences of $\mathbb{R}^{|A|\times n_{L+1}}$ valued random variables (and we 
shall often use the indices $\alpha h$ to mean $(\alpha,h)\in A\times\{1,\ldots,n_{L+1}\}$):
\begin{itemize}
\item The large deviation principle of the sequence
\begin{equation}\label{eq:ZLDP}
\big\{\big(Z_h^{(L+1)}(x_\alpha)/\sqrt{v_n^*}\big)_{\alpha h}\big\}_n,
\end{equation}
with speed $v_n^*$ 
(see Theorem \ref{th:finitedimensional-LD}).
\item For every sequence of positive numbers $\{a_n\}_n$ such that
\begin{equation}\label{eq:MD-conditions}
		a_n\to 0\quad\mbox{and}\quad a_nv_n^*\to\infty,\quad\text{as $n\to\infty$},
\end{equation}
the large deviation principle of the sequence
\begin{equation}\label{eq:ZMDP}
\big\{\big(\sqrt{a_n}Z_h^{(L+1)}(x_\alpha)\big)_{\alpha h}\big\}_n, 
\end{equation}
with speed $1/a_n$
(see Theorem \ref{th:finitedimensional-MD}).
\end{itemize}
We remark that the sequences in \eqref{eq:ZLDP} and \eqref{eq:ZMDP} converge to the null matrix 
$0\in\mathbb{R}^{|A|\times n_{L+1}}$, where $|A|$ denotes the cardinality of $A$. The class of large deviation principles
in the second result is known in literature as a moderate deviation principle and fills the gap between the convergence 
in probability to $0\in\mathbb{R}^{|A|\times n_{L+1}}$ of the sequence in \eqref{eq:ZLDP},
 and the convergence in law of the sequence $\{(Z_h^{(L+1)}(x_\alpha))_{\alpha h}\}_n$
to a Gaussian vector. We remark that moderate deviation estimates (and concentration inequalities) for the 
output of a Gaussian neural network might be proved even applying the classical theory on the subject developed in 
\cite{SaulisStatulevicius}
. To apply this theory, one needs a fine study of the cumulants of the output.
To this aim, the main findings in \cite{Hanin-unpublished} seem to be not useful as they do not provide the constants involved 
in the big $O$ functions which give the rate of the cumulants. This could be a topic for a future research.

The proofs of our main results proceed by induction on the number of layers and combine a representation of Gaussian fully 
connected deep neural networks (see Lemma \ref{lem:equality-in-law}) with a large deviation principle on product spaces proved 
in \cite{Chaganty} (see Proposition \ref{prop:Chaganty}). These results hold when the pre-activation function is bounded and 
continuous,  and therefore exclude some important classes of pre-activations such as the ReLU function. However, we are able to 
provide large and moderate deviation principles for Gaussian fully connected deep neural networks with a single input and a 
ReLU pre-activation (see Section \ref{sec:ReLUoneinput}). Indeed, if $|A|=1$, some technical difficulties encountered in the 
general case can be overcome by means of ad hoc arguments. Recently our result for the case
of ReLU function has been generalized in \cite{Vogel}, where the author considers Gaussian fully connected deep neural networks 
with linearly growing pre-activation functions. 

For the case of Gaussian shallow neural networks, i.e., if $L=1$, the model is much more simple.
In such a case we are able to prove large and moderate deviations for the output and its sensitivity (i.e.,  the derivative of 
the output with respect to the input), under quite general assumptions on $\sigma$ (see Propositions \ref{prop:derivative-LD}, 
\ref{prop:derivative-MD} and \ref{prop:for-the-finiteness-around-the-origin}).


The paper is structured as follows. In Section \ref{sec:main-results-and-remarks} we give the statements of the 
main results, together with the preliminary notation and some remarks. The proofs of the main results are presented in Section \ref{sec:proofs-main-results}. The particular case of a Gaussian fully connected deep neural network with a single input and 
ReLU pre-activation function is treated in Section \ref{sec:ReLUoneinput}. Finally, large and moderate deviations for Gaussian 
shallow neural networks and their sensitivities are presented in Section \ref{sec:L=1}.

\section{Main results: statements and remarks}\label{sec:main-results-and-remarks}

We start with some preliminary notation and the hypotheses.
In our main results we assume that the following Condition \ref{cond:on-sigma} holds.

\begin{condition}\label{cond:on-sigma}
	The function $\sigma$ is continuous and bounded. In particular we set
	$$\|\sigma\|_\infty:=\sup_{x\in\mathbb{R}}|\sigma(x)|.$$
\end{condition}

As already mentioned in the introduction, we take
$$T:=\{x_\alpha\}_{\alpha\in A},\quad \mbox{with}\ x_\alpha=(x_{\alpha,1},\ldots,x_{\alpha,n_0})\in\mathbb{R}^{n_0},$$
for some finite set $A$. 
Moreover, we shall use the simplified notation $\underline{x}$ to mean $T=\{x_\alpha\}_{\alpha\in A}$.


It is well-known (see, e.g., Section I.2 in \cite{Bhatia}) that, for every symmetric 
positive semidefinite matrix $q=(q_{\alpha\beta})_{\alpha,\beta\in A}$,  there exists a unique
symmetric positive semidefinite matrix $q^\#=(q_{\alpha\beta}^\#)_{\alpha,\beta\in A}$ such that 
$$q^\# q^\#=q,\quad \mbox{i.e.}\ q_{\alpha\beta}=\sum_{\gamma\in A}q_{\alpha\gamma}^\# q_{\gamma\beta}^\#\ (\mbox{for all}\ \alpha,\beta\in A).$$
Throughout this paper we use the notation $\mathbf{1}\in\mathbb{R}^{|A|\times|A|}$ for the matrix with all entries equal to 1.
Moreover we denote by $\mathcal{S}_{|A|,C_b}$ be the family of symmetric positive semidefinite matrices 
$q=(q_{\alpha\beta})_{\alpha,\beta\in A}$ such that 
$$q-C_b\mathbf{1}=(q_{\alpha\beta}-C_b)_{\alpha,\beta\in A}$$
is again a symmetric positive semidefinite matrix. 

Let $(N_\gamma)_{\gamma\in A}$,  $|A|<\infty$,  be a family of independent standard Normal distributed random variables.  Hereafter we consider the function $\kappa(\cdot;q)$ (where $q=(q_{\alpha\beta})_{\alpha,\beta\in A}\in\mathbb{R}^{|A|\times|A|}$ is a symmetric and 
positive semidefinite matrix) defined by
\begin{equation}\label{eq:def-kappa}
	\kappa(\eta;q):=\log\mathbb{E}\big[\exp\big(\sum_{\alpha,\beta\in A}\eta_{\alpha\beta}
	\sigma(\langle q_{\alpha\cdot}^\#,N_\cdot\rangle)\sigma(\langle q_{\beta\cdot}^\#,N_\cdot\rangle)\big)\big],
\mbox{for every}\ \eta=(\eta_{\alpha\beta})_{\alpha,\beta\in A}\in\mathbb{R}^{|A|\times|A|},
\end{equation}
where
$\langle q_{\gamma\cdot}^\#,N_\cdot\rangle:=\sum_{\gamma^\prime\in A}q_{\gamma\gamma^\prime}^\# N_{\gamma^\prime}$ (for every $\gamma\in A$).
We note that, under the Condition \ref{cond:on-sigma}, $\kappa(\cdot;q)$ assumes finite values;
moreover, it is easy to check that $(\eta,q)\mapsto\kappa(\eta;q)$ is continuous, indeed $q\mapsto q^\#$ is continuous (see, e.g., 
Theorem X.1.1 in \cite{Bhatia}).

In what follows we also consider the Fenchel-Legendre transform of $\kappa$, i.e. the function $\kappa^*(\cdot;q)$ defined by
\begin{equation}\label{eq:def-kappa-star}
	\kappa^*(y;q)=\sup_{\eta\in\mathbb{R}^{|A|\times|A|}}\{\langle\eta,y\rangle-\kappa(\eta;q)\},
	\quad\mbox{where}\ \langle\eta,y\rangle=\sum_{\alpha,\beta\in A}\eta_{\alpha\beta}y_{\alpha\beta};
\end{equation}
then we have $\kappa^*(y;q)=0$ if and only if $y=y(q)$, where $y(q):=(y_{\alpha\beta}(q))_{\alpha,\beta\in A}$ and
\begin{equation}\label{eq:unique-zero-LT}
	y_{\alpha\beta}(q):=\left.\frac{\partial\kappa(\eta;q)}{\partial\eta_{\alpha\beta}}\right|_{\eta=0}
	=\mathbb{E}\big[\sigma(\langle q_{\alpha\cdot}^\#,N_\cdot\rangle)\sigma(\langle q_{\beta\cdot}^\#,N_\cdot\rangle)\big]
\end{equation}
(here $0$ in the null matrix in $\mathbb{R}^{|A|\times|A|}$).

Throughout this paper we assume the following condition on the widths $n_\ell=n_\ell(n)$ ($\ell=1,\ldots,L$), which we recall tend to 
infinity as $n\to\infty$.

\begin{condition}\label{cond:on-n1-nL}
	There exists $\widehat{\ell}\in\{1,\ldots,L\}$ such that, for some $\gamma_1,\ldots,\gamma_L\in[1,\infty]$,
	$$\lim_{n\to\infty}\frac{n_\ell(n)}{n_{\widehat{\ell}}(n)}=\gamma_\ell\quad(\mbox{for all}\ \ell=1,\ldots,L)$$
    In what follows we simply write $v_n^*$ in place of $n_{\widehat{\ell}}(n)$.
\end{condition}

Note that Condition \ref{cond:on-n1-nL} always holds when $L=1$ (we have $\widehat{\ell}=1$ and $\gamma_1=1$); moreover, if Condition 
\ref{cond:on-n1-nL} holds, we have $\gamma_{\widehat{\ell}}=1$.

We conclude with some further notation. Let
$g^{(0)}(\underline{x})=(g_{\alpha\beta}^{(0)}(x_\alpha,x_\beta))_{\alpha,\beta\in A}\in\mathcal{S}_{|A|,C_b}$ defined by
\begin{equation}\label{eq:initial-variance}
	g_{\alpha\beta}^{(0)}(x_\alpha,x_\beta):=\mathrm{Cov}(Z_i^{(1)}(x_\alpha),Z_i^{(1)}(x_\beta))
	=C_b+\frac{C_W}{n_0}\sum_{r=1}^{n_0}x_{\alpha,r}x_{\beta,r}
\end{equation}
(indeed the covariance in \eqref{eq:initial-variance} does not depend on $i\in\{1,\ldots,n_1\}$). Moreover, let
$\widehat{g}_{\underline{x}}^{(L)}\in\mathcal{S}_{|A|,C_b}$ be defined by recurrence 
(on $L\geq 1$) by
\begin{equation}\label{eq:Hanin-covariance}
	\widehat{g}_{\underline{x}}^{(\ell)}=C_b\mathbf{1}+C_Wy(\widehat{g}_{\underline{x}}^{(\ell-1)})\quad(\mbox{for all}\ \ell=1,\ldots,L),
\end{equation}
where $y(q)$ is defined by \eqref{eq:unique-zero-LT} and $\widehat{g}_{\underline{x}}^{(0)}=g^{(0)}(\underline{x})$ as in 
\eqref{eq:initial-variance}.

Now we are ready to present the statements of the main theorems. All the preliminaries on large deviations will be given 
in Section \ref{sec:proofs-main-results}. We shall use the acronymous LDP to mean \emph{large deviation principle}.

\begin{theorem}[Large deviations]\label{th:finitedimensional-LD}
	Assume that Conditions \ref{cond:on-sigma} and \ref{cond:on-n1-nL} hold. Then the sequence \\
	$\{(Z_h^{(L+1)}(x_\alpha)/\sqrt{v_n^*})_{\alpha h}\}_n$ satisfies the 
	LDP on $\mathbb{R}^{|A|\times n_{L+1}}$, with speed $v_n^*$ and good rate function $I_{Z^{(L+1)}(\underline{x})}$ defined by
	\begin{equation}\label{eq:finitedimensional-rf}
		I_{Z^{(L+1)}(\underline{x})}(z):=\inf_{g^{(L)}\in\mathcal{S}_{|A|,C_b},r\in\mathbb{R}^{|A|\times n_{L+1}}}
		\{I_{G^{(L)}(\underline{x})}(g^{(L)})+\|r\|^2/2:g^{(L),\#}r=z\},
	\end{equation}
	where: $\|\cdot\|$ is the Euclidean norm in $\mathbb{R}^{|A|\times n_{L+1}}$, $I_{G^{(L)}(\underline{x})}$ is defined by
	\begin{equation}\label{eq:finitedimensional-rf-auxiliary}
		I_{G^{(L)}(\underline{x})}(g^{(L)}):=\inf\big\{\sum_{\ell=1}^LJ(g^{(\ell)}|g^{(\ell-1)}):
		g^{(0)}=g^{(0)}(\underline{x}),g^{(1)},\ldots,g^{(L-1)}\in\mathcal{S}_{|A|,C_b}\big\}
	\end{equation}
	(for $g^{(L)}\in\mathcal{S}_{|A|,C_b}$), $g^{(0)}(\underline{x})$ is defined by \eqref{eq:initial-variance}, 
	$J(\cdot|\cdot)$ is defined by
	\begin{equation}\label{eq:finitedimensional-rf-auxiliary-bis}
		J(g^{(\ell)}|g^{(\ell-1)}):=\left\{\begin{array}{ll}
			\gamma_\ell\kappa^*(\frac{g^{(\ell)}-C_b\mathbf{1}}{C_W};g^{(\ell-1)})&\ \mbox{if}\ \gamma_\ell<\infty\\
			\Delta(g^{(\ell)};C_b\mathbf{1}+C_Wy(g^{(\ell-1)}))&\ \mbox{if}\ \gamma_\ell=\infty
		\end{array}\right.
	\end{equation}
	(for $g^{(\ell-1)},g^{(\ell)}\in\mathcal{S}_{|A|,C_b}$), $\kappa^*(\cdot,\cdot)$ is defined by
	\eqref{eq:def-kappa-star}, $\Delta(\cdot,\cdot)$ is defined just before Lemma \ref{lem:two-speed-functions} and $y(g^{(\ell-1)})=(y_{\alpha\beta}(g^{(\ell-1)}))_{\alpha,\beta\in A}$ is defined by \eqref{eq:unique-zero-LT}.
\end{theorem}

\begin{theorem}[Moderate deviations]\label{th:finitedimensional-MD}
	Assume that Conditions \ref{cond:on-sigma} and \ref{cond:on-n1-nL} hold. Then, for every sequence of positive numbers 
	$\{a_n\}_n$ such that \eqref{eq:MD-conditions} holds, the sequence
	$\{(\sqrt{a_n}Z_h^{(L+1)}(x_\alpha))_{\alpha h}\}_n$ satisfies the LDP 
	on $\mathbb{R}^{|A|\times n_{L+1}}$, with speed $1/a_n$ and good rate function $\widetilde{I}_{Z^{(L+1)}(\underline{x})}$
	defined by
	\begin{equation}\label{eq:finitedimensional-rf-MD}
		\widetilde{I}_{Z^{(L+1)}(\underline{x})}(z):=\inf_{r\in\mathbb{R}^{|A|\times n_{L+1}}}
		\{\|r\|^2/2:\widehat{g}_{\underline{x}}^{(L),\#}r=z\},
	\end{equation}
	where $\|\cdot\|$ is the Euclidean norm in $\mathbb{R}^{|A|\times n_{L+1}}$ and $\widehat{g}_{\underline{x}}^{(L)}$ is 
	defined by recurrence (on $L\geq 1$) by \eqref{eq:Hanin-covariance}, with 
	$\widehat{g}_{\underline{x}}^{(0)}=g^{(0)}(\underline{x})$ as in \eqref{eq:initial-variance}. Thus, if
	$\widehat{g}_{\underline{x}}^{(L)}$ is invertible (and therefore $\widehat{g}_{\underline{x}}^{(L),\#}$ is also invertible),
	we have
	$$\widetilde{I}_{Z^{(L+1)}(\underline{x})}(z)=\|(\widehat{g}_{\underline{x}}^{(L),\#})^{-1}z\|^2/2.$$
\end{theorem}

Theorem \ref{th:finitedimensional-MD} provides a class of LDPs which fills the gap between the convergence to zero governed 
by the LDP in Theorem \ref{th:finitedimensional-LD}, and the weak convergence in Theorem 1.2 of
\cite{Hanin} cited above; these two asymptotic regimes correspond to the cases $a_n=1/v_n^*$ and $a_n=1$, respectively (in both cases 
one condition in \eqref{eq:MD-conditions} holds, and the other one fails).

We conclude this section with some remarks.

\begin{remark}\label{rem:finitedimensional-rf-LD-more-explicit}
	We can say that $I_{Z^{(L+1)}(\underline{x})}(z)=0$ if and only if $z=0\in\mathbb{R}^{|A|\times n_{L+1}}$ because we have
	$$I_{G^{(L)}(\underline{x})}(g^{(L)})+\|r\|^2/2=0$$ 
	when $r=0\in\mathbb{R}^{|A|\times n_{L+1}}$ and $g^{(L)}=\widehat{g}_{\underline{x}}^{(L)}$, where $\widehat{g}_{\underline{x}}^{(L)}\in\mathcal{S}_{|A|,C_b}$ is defined as in the statement of Theorem 
	\ref{th:finitedimensional-MD}.
	Indeed, for every $\ell=1,\ldots,L$, one can check with some computations that 
	$J(\widehat{g}_{\underline{x}}^{(\ell)}|\widehat{g}_{\underline{x}}^{(\ell-1)})=0$ and, by 
	\eqref{eq:finitedimensional-rf-auxiliary}, we have $I_{G^{(L)}(\underline{x})}(\widehat{g}_{\underline{x}}^{(L)})=0$.
	We can also say that $\widetilde{I}_{Z^{(L+1)}(\underline{x})}(z)=0$ if and only if $z=0\in\mathbb{R}^{|A|\times n_{L+1}}$.
\end{remark}

\begin{remark}\label{rem:link-with-Hanin-covariance}
	The matrix $\widehat{g}_{\underline{x}}^{(\ell)}$ in \eqref{eq:Hanin-covariance} coincides with 
	$K^{(\ell+1)}=(K_{\alpha\beta}^{(\ell+1)})_{\alpha,\beta\in A}$ in \cite{Hanin} (see eqs. (1.7) and (1.8) in that reference).
\end{remark}

\begin{remark}\label{rem:more-explicit-|A|=1}
	If $|A|=1$ we simply have $\mathcal{S}_{|A|,C_b}=[C_b,\infty)$ and $x\in\mathbb{R}^{n_0}$ in place of $\underline{x}$. 
	Then, if we specialize \eqref{eq:finitedimensional-rf} to this case, we get
	$$I_{Z^{(L+1)}(x)}(z):=\inf_{g^{(L)}\geq C_b}\{I_{G^{(L)}(x)}(g^{(L)})+\|z\|^2/(2g^{(L)})\}
	\quad (z\in\mathbb{R}^{n_{L+1}}).$$
	Note that we can have $g^{(L)}=0$ if $C_b=0$ and, if we consider the rule $\frac{0}{0}=0$ as usual, the 
	argument of the infimum above computed at $g^{(L)}=0$ is	
	$$\left.\{I_{G^{(L)}(x)}(g^{(L)})+\|z\|^2/(2g^{(L)})\}\right|_{g^{(L)}=0}=\left\{\begin{array}{ll}
		I_{G^{(L)}(x)}(0)&\ \mbox{if}\ z=0\\
		\infty&\ \mbox{if}\ z\neq 0.
	\end{array}\right.$$
	Moreover we have $\widehat{g}_x^{(L)}\geq C_b\geq 0$, and \eqref{eq:finitedimensional-rf-MD} yields
    $$\widetilde{I}_{Z^{(L+1)}(x)}(z):=\|z\|^2/(2\widehat{g}_x^{(L)})\quad (z\in\mathbb{R}^{n_{L+1}}).$$
    Finally we also remark that, if $|A|=1$ and $\widehat{g}_x^{(L)}=0\in\mathbb{R}$ (this can happen if $C_b=0$), we have
    $\widetilde{I}_{Z^{(L+1)}(x)}(z)=0$ if $z=0\in\mathbb{R}^{n_{L+1}}$, and $\widetilde{I}_{Z^{(L+1)}(x)}(z)=\infty$
    otherwise.
\end{remark}

\section{Proofs of the main results}\label{sec:proofs-main-results}
In this section we recall some preliminaries on large deviations (basic definitions and results), we provide
an important representation lemma and, finally, we present the proofs of Theorems \ref{th:finitedimensional-LD} 
and \ref{th:finitedimensional-MD}.

\subsection{Preliminaries on large deviations}
We start with the basic definition of large deviation principle (see e.g. \cite{DemboZeitouni}), and other related concepts.
In view of what follows we present these definitions by referring to sequences of probability measures.

\begin{definition}
	A sequence of positive numbers $\{v_n:n\geq 1\}$ such that $v_n\to\infty$ (as $n\to\infty$) is 
	called \emph{speed function}, and a lower semicontinuous function $I:\mathcal{X}\to[0,\infty]$ is called 
	\emph{rate function}. Let $\mathcal{X}$ be a topological space, and let $\{\pi_n\}_n$ be a sequence of probability measures 
	on $\mathcal{X}$ (equipped with its completed Borel $\sigma$-field). Then $\{\pi_n\}_n$ satisfies the \emph{large deviation 
	principle} (LDP from now on) 
	on $\mathcal{X}$, with speed $v_n$ and rate function $I$ if
	$$\limsup_{n\to\infty}\frac{1}{v_n}\log\pi_n(C)\leq-\inf_{x\in C}I(x)\quad\mbox{for all closed sets}\ C,$$
	and
	$$\liminf_{n\to\infty}\frac{1}{v_n}\log\pi_n(O)\geq-\inf_{x\in O}I(x)\quad\mbox{for all open sets}\ O.$$
	Moreover the rate function $I$ is said to be \emph{good} if, for every $\lambda\geq 0$, the level set 
	$\{x\in\mathbb{R}:I(x)\leq\lambda\}$ is compact. If the upper bound above holds for compact sets only, then we say that
	$\{\pi_n\}_n$ satisfies the \emph{weak large deviation principle} (WLDP from now on) on 
	$\mathcal{X}$.
\end{definition}

The results in this paper are stated for sequences of random variables $\{R_n\}_n$, say, defined on the same probability 
space $(\Omega,\mathcal{F},P)$, and taking values on the same topological space $\mathcal{X}$. Indeed we refer to the sequence
of laws, i.e. to the framework of the above definition with $\pi_n=P(R_n\in\cdot)$.

Throughout this paper we refer to some well-known large deviation results: the \emph{G\"artner Ellis Theorem} 
(see e.g. Theorem 2.3.6(c) in \cite{DemboZeitouni}), and the \emph{contraction principle} (see e.g. Theorem 4.2.1 in 
\cite{DemboZeitouni}). We also refer to the concept of \emph{essentially smooth} function (see e.g. Definition 2.3.5 in 
\cite{DemboZeitouni}).

An important result used in this paper is Theorem 2.3 in \cite{Chaganty}. We recall its statement.

\begin{proposition}\label{prop:Chaganty}
	Let $\Omega_1$ and $\Omega_2$ be two Polish spaces. Let $\{\pi_n:n\geq 1\}$ be a sequence of probability measures on 
	$\mathcal{X}=\Omega_1\times\Omega_2$, and denote the sequences of marginal distributions by 
	$\{\pi_{1,n}:n\geq 1\}$ on $\Omega_1$ and $\{\pi_{2,n}:n\geq 1\}$ on $\Omega_2$. We assume that the following conditions hold.
	\begin{enumerate}
		\item The sequence $\{\pi_{1,n}:n\geq 1\}$ satisfies the LDP on $\Omega_1$, with speed $v_n$ and a good rate function 
		$I_1$.
		\item If $\{x_{1,n}:n\geq 1\}\subset\Omega_1$ and $x_{1,n}\to x_1\in\Omega_1$, then the sequence of conditional distributions
		$\{\pi_n^{2|1}(dx_2|x_{1,n}):n\geq 1\}$ satisfies the LDP on $\Omega_2$, with speed $v_n$ and rate function $J(\cdot|x_1)$.
		\item $(x_1,x_2)\mapsto J(x_2|x_1)$ is lower semicontinuous.
	\end{enumerate}
	Then $\{\pi_n:n\geq 1\}$ satisfies the WLDP on $\Omega_1\times\Omega_2$, with speed $v_n$ and rate function $I$ defined by
	$$I(x_1,x_2):=J(x_2|x_1)+I_1(x_1).$$
	Moreover:
	$\{\pi_{2,n}:n\geq 1\}$ satisfies the LDP on $\Omega_2$, with speed $v_n$ and rate function $I_2$ defined by
	$$I_2(x_2):=\inf_{x_1\in\Omega_1}\{I(x_1,x_2)\}=\inf_{x_1\in\Omega_1}\{J(x_2|x_1)+I_1(x_1)\};$$
	$\{\pi_n:n\geq 1\}$ satisfies the LDP on $\Omega_1\times\Omega_2$ if the rate function $I$ is good and, in such a case, 
	the rate function $I_2$ is also good.
\end{proposition}

We also recall another useful related result (for its proof see Lemma 2.6 in \cite{Chaganty}).

\begin{lemma}\label{lem:Chaganty}
	Consider the same hypotheses and notation of Proposition \ref{prop:Chaganty}. Assume that, for every $a\geq 0$ and for 
	every compact subset $K_1$ of $\Omega_1$, the set
	$$\bigcup_{x_1\in K_1}\{x_2\in\Omega_2:J(x_2|x_1)\leq a\}$$
	is a compact subset of $\Omega_2$. Then the rate function $I$ in Proposition \ref{prop:Chaganty} is good.
\end{lemma}

We conclude this section with a useful lemma. We consider a Polish space $\Pi$ and, for any given $\bar{r}\in\Pi$, let 
$\Delta(\cdot;\bar{r}):\Pi\to[0,\infty]$ be the function defined by
$$\Delta(r;\bar{r}):=\left\{\begin{array}{ll}
	0&\ \mbox{if}\ r=\bar{r}\\
	\infty&\ \mbox{if}\ r\neq\bar{r}.
\end{array}\right.$$
In particular we remark that $\Delta(\cdot;\bar{r})$ is trivially a good rate function; indeed every level set is compact because
$$\{r\in\mathbb{R}:\Delta(r;\bar{r})\leq a\}=\{\bar{r}\}\quad(\mbox{for all}\ a\geq 0).$$
In what follows we consider this function with $\Pi=\mathcal{S}_{|A|,C_b}\subset\mathbb{R}^{|A|\times|A|}$.

\begin{lemma}[Lemma 1 in \cite{GiulianoMacciPacchiarotti}]\label{lem:two-speed-functions}
	Let $\{\psi_n\}_n$ be a sequence of probability measures on some Polish space $\Pi$ that satisfies the LDP with speed $s_n$ and good rate 
	function $H$, which uniquely vanishes at some $r_0$. Moreover let $t_n$ be another speed function such that $\frac{s_n}{t_n}\to\infty$. 
	Then $\{\psi_n\}_n$ satisfies the LDP with speed $t_n$ and good rate function $\Delta(\cdot;r_0)$.
\end{lemma}

\subsection{An important representation lemma}
The proofs of the main results are based on a different representation of the model. Such a representation is based on a recursive
approach described without any mathematical details in Section 1.4 in \cite{Hanin} (see the part with eqs. (1.9) and (1.10) in that 
reference). Roughly speaking one can say that the output of the network
(layer $L+1$) is centered multivariate Normal distributed, with random covariance matrix which depends on the random variables 
involved in the previous layers.

The aim of this section is to present a lemma with the mathematical details of the recursive approach described in \cite{Hanin}. 
This representation holds for some fixed choices of $n_1,\ldots,n_L$; however, for the future use of this lemma, we write 
$n_1(n),\ldots,n_L(n)$ in place of $n_1,\ldots,n_L$ as in the most of the paper.

\begin{lemma}\label{lem:equality-in-law}
	Let $\{N_{j\alpha}^{(\ell)}:j,\ell\geq 1,\alpha\in A\}$ be a family of independent standard Normal distributed random variables,
	and let $\{G_n^{(\ell)}(\underline{x}):n,\ell\geq 1\}$ be the $\mathcal{S}_{|A|,C_b}$-valued random variables defined (by recurrence on 
	$\ell\geq 1$) by $G_n^{(\ell)}(\underline{x}):=(G_{n;\alpha\beta}^{(\ell)}(\underline{x}))_{\alpha,\beta\in A}$, where
	\begin{equation}\label{eq:random-variances}
		G_{n;\alpha\beta}^{(\ell)}(\underline{x}):=
		C_b+\frac{C_W}{n_\ell(n)}\sum_{j=1}^{n_\ell(n)}
		\sigma(\langle G_{n;\alpha\cdot}^{(\ell-1),\#}(\underline{x}),N_{j\cdot}^{(\ell)}\rangle)
		\sigma(\langle G_{n;\beta\cdot}^{(\ell-1),\#}(\underline{x}),N_{j\cdot}^{(\ell)}\rangle),
	\end{equation}
		$\langle G_{n;\gamma\cdot}^{(\ell-1),\#}(\underline{x}),N_{j\cdot}^{(\ell)}\rangle:=
		\sum_{\gamma^\prime\in A}G_{n;\gamma\gamma^\prime}^{(\ell-1),\#}(\underline{x}) N_{j\gamma^\prime}^{(\ell)}$ (for every $\gamma\in A$),
	and $G_n^{(0)}(\underline{x})=g^{(0)}(\underline{x})$ as in \eqref{eq:initial-variance}. Then, for $L\geq 1$,
	$$(Z_h^{(L+1)}(x_\alpha))_{\alpha h}\stackrel{\mathrm{law}}{=}
	\big(\sum_{\gamma\in A}G_{n;\alpha\gamma}^{(L),\#}(\underline{x})N_{h\gamma}^{(L+1)}\big)_{\alpha h}.$$
\end{lemma}

\begin{proof}
	Throughout this appendix we simply write $n_\ell$ in place of $n_\ell(n)$. Moreover, we can say that $\{G_n^{(\ell)}(\underline{x}):n,\ell\geq 1\}$
	are $\mathcal{S}_{|A|,C_b}$-valued random variables by construction. In what follows we prove the equality in law by induction 
	on $L$, showing that the moment generating functions (with argument $\theta\in\mathbb{R}^{|A|\times n_{L+1}}$) coincide.
	
	We start with the case $L=0$, i.e.,
	$$\Big(\Big(Z_h^{(1)}(x_\alpha)\Big)_{\alpha\in A}\Big)_{h=1,\ldots,n_1}\stackrel{\mathrm{law}}{=}
	\Big(\Big(\sum_{\gamma\in A}G_{n;\alpha\gamma}^{(0),\#}(\underline{x})N_{h\gamma}^{(1)}\Big)_{\alpha\in A}\Big)_{h=1,\ldots,n_1},$$
	where $G_n^{(0)}(\underline{x})=g^{(0)}(\underline{x})$ as in \eqref{eq:initial-variance}, and therefore 
	$G_n^{(0),\#}(\underline{x})=g^{(0),\#}(\underline{x})$. By \eqref{eq:Hanin-model} and some manipulations we have
	\begin{multline*}
		\mathbb{E}\Big[\exp\Big(\sum_{\alpha\in A}\sum_{h=1}^{n_1}\theta_{\alpha h}Z_h^{(1)}(x_\alpha)\Big)\Big]
		\\
		=\exp\Big(\frac{1}{2}\Big(C_b\sum_{h=1}^{n_1}\sum_{\alpha,\beta\in A}\theta_{\alpha h}\theta_{\beta h}
		+\frac{C_W}{n_0}\sum_{h=1}^{n_1}\sum_{r=1}^{n_0}\sum_{\alpha,\beta\in A}\theta_{\alpha h}\theta_{\beta h}x_{\alpha,r}x_{\beta,r}\Big)\Big).
	\end{multline*}
	On the other hand we also have (in the final equality we take into account \eqref{eq:initial-variance})
	\begin{multline*}
		\mathbb{E}\Big[\exp\Big(\sum_{\alpha\in A}\sum_{h=1}^{n_1}\theta_{\alpha h}
		\sum_{\gamma\in A}G_{n;\alpha\gamma}^{(0),\#}(\underline{x})N_{h\gamma}^{(1)}\Big)\Big]
		=\exp\Big(\frac{1}{2}\sum_{h=1}^{n_1}\sum_{\gamma\in A}\Big(\sum_{\alpha\in A}\theta_{\alpha h}g_{\alpha\gamma}^{(0),\#}(\underline{x})\Big)^2\Big)
		\\
		=\exp\Big(\frac{1}{2}\sum_{h=1}^{n_1}\sum_{\alpha,\beta\in A}\theta_{\alpha h}\theta_{\beta h}g_{\alpha\beta}^{(0)}(\underline{x})\Big)
		=\exp\Big(\frac{1}{2}\sum_{h=1}^{n_1}\sum_{\alpha,\beta\in A}\theta_{\alpha h}\theta_{\beta h}
		\Big(C_b+\frac{C_W}{n_0}\sum_{r=1}^{n_0}x_{\alpha,r}x_{\beta,r}\Big)\Big).
	\end{multline*}
	So the case $L=0$ is proved.
	
	Now we assume that the statement for $L$ is proved, i.e.,
	$$\Big(\Big(Z_h^{(L)}(x_\alpha)\Big)_{\alpha\in A}\Big)_{h=1,\ldots,n_L}\stackrel{\mathrm{law}}{=}
	\Big(\Big(\sum_{\gamma\in A}G_{n;\alpha\gamma}^{(L-1),\#}(\underline{x})N_{h\gamma}^{(L)}\Big)_{\alpha\in A}\Big)_{h=1,\ldots,n_L},$$
	and we prove the statement for $L+1$, i.e.,
	$$\Big(\Big(Z_h^{(L+1)}(x_\alpha)\Big)_{\alpha\in A}\Big)_{h=1,\ldots,n_{L+1}}\stackrel{\mathrm{law}}{=}
	\Big(\Big(\sum_{\gamma\in A}G_{n;\alpha\gamma}^{(L),\#}(\underline{x})N_{h\gamma}^{(L+1)}\Big)_{\alpha\in A}\Big)_{h=1,\ldots,n_{L+1}}.$$
	By \eqref{eq:Hanin-model} and some manipulations 
	we have
	\begin{multline*}
		\mathbb{E}\Big[\exp\Big(\sum_{\alpha\in A}\sum_{h=1}^{n_{L+1}}\theta_{\alpha h}Z_h^{(L+1)}(x_\alpha)\Big)\Big]
		\\
		=\mathbb{E}\Big[\exp\Big(\frac{C_b}{2}\sum_{h=1}^{n_{L+1}}\Big(\sum_{\alpha\in A}\theta_{\alpha h}\Big)^2\Big)
		\exp\Big(\frac{C_W}{2n_L}\sum_{h=1}^{n_{L+1}}\sum_{j=1}^{n_L}\Big(\sum_{\alpha\in A}\theta_{\alpha h}\sigma(Z_j^{(L)}(x_\alpha))\Big)^2\Big)\Big]\\
		=\exp\Big(\frac{C_b}{2}\sum_{h=1}^{n_{L+1}}\sum_{\alpha,\beta\in A}\theta_{\alpha h}\theta_{\beta h}\Big)
		\mathbb{E}\Big[\exp\Big(\frac{C_W}{2n_L}\sum_{h=1}^{n_{L+1}}\sum_{j=1}^{n_L}
		\sum_{\alpha,\beta\in A}\theta_{\alpha h}\theta_{\beta h}\sigma(Z_j^{(L)}(x_\alpha))\sigma(Z_j^{(L)}(x_\beta))\Big)\Big];
	\end{multline*}
	thus, by some manipulations (in the first equality we take into account the inductive hypothesis, in the second equality we take into account 
	by \eqref{eq:random-variances}), we obtain
	\begin{multline*}
		\mathbb{E}\Big[\exp\Big(\sum_{\alpha\in A}\sum_{h=1}^{n_{L+1}}\theta_{\alpha h}Z_h^{(L+1)}(x_\alpha)\Big)\Big]\\
		=\mathbb{E}\Big[\exp\Big(\frac{1}{2}\sum_{h=1}^{n_{L+1}}\sum_{\alpha,\beta\in A}\theta_{\alpha h}\theta_{\beta h}
		\Big(C_b+\frac{C_W}{n_L}\sum_{j=1}^{n_L}\sigma\Big(\sum_{\gamma\in A}G_{n;\alpha\gamma}^{(L-1),\#}(\underline{x})
		N_{j\gamma}^{(L)}\Big)\Big.\Big.\Big.\\
		\Big.\Big.\Big.\times\sigma\Big(\sum_{\gamma\in A}G_{n;\beta\gamma}^{(L-1),\#}(\underline{x})N_{j\gamma}^{(L)}\Big)\Big)\Big)\Big]
		=\mathbb{E}\Big[\exp\Big(\frac{1}{2}\sum_{h=1}^{n_{L+1}}\sum_{\alpha,\beta\in A}\theta_{\alpha h}\theta_{\beta h}
		G_{n;\alpha\beta}^{(L)}(\underline{x})\Big)\Big].
	\end{multline*}
	On the other hand we also have 
	\begin{multline*}
		\mathbb{E}\Big[\exp\Big(\sum_{\alpha\in A}\sum_{h=1}^{n_{L+1}}\theta_{\alpha h}
		\sum_{\gamma\in A}G_{n;\alpha\gamma}^{(L),\#}(\underline{x})N_{h\gamma}^{(L+1)}\Big)\Big]
		\\
		=\mathbb{E}\Big[\exp\Big(\frac{1}{2}\sum_{h=1}^{n_{L+1}}\sum_{\gamma\in A}\Big(\sum_{\alpha\in A}\theta_{\alpha h}
		G_{n;\alpha\gamma}^{(L),\#}(\underline{x})\Big)^2\Big)\Big]
		=\mathbb{E}\Big[\!\exp\Big(\frac{1}{2}\sum_{h=1}^{n_{L+1}}\sum_{\alpha,\beta\in A}\theta_{\alpha h}\theta_{\beta h}
		G_{n;\alpha\beta}^{(L)}(\underline{x})\Big)\!\Big].
	\end{multline*}
	The statement for $L+1$ is proved.
\end{proof}

\subsection{Proof of Theorem \ref{th:finitedimensional-LD}}
	We prove the theorem by induction on $L$. Moreover, by Lemma \ref{lem:equality-in-law}, we refer to the random variables 
	$\big(\sum_{\gamma\in A}G_{n;\alpha\gamma}^{(L),\#}(\underline{x})N_{h\gamma}^{(L+1)}\big)_{\alpha h}$
	in place of the random variables $(Z_h^{(L+1)}(x_\alpha))_{\alpha h}$.
	
	Actually we prove by induction on $L$ that
	$$(\bullet):\ \left\{\begin{array}{l}
			\mbox{$\{G_n^{(L)}(\underline{x})\}_n$ satisfies the LDP on $\mathcal{S}_{|A|,C_b}$,}\\
			\mbox{with speed $v_n^*$ and good rate function $I_{G^{(L)}(\underline{x})}$ in \eqref{eq:finitedimensional-rf-auxiliary}.}
		\end{array}\right.$$
	Indeed, if $(\bullet)$ holds, we can conclude as follows: 
	$\{(N_{\alpha h}^{(L+1)}/\sqrt{v_n^*})_{\alpha h}\}_n$ satisfies the LDP (on\\
	$\mathbb{R}^{|A|\times n_{L+1}}$)
	with good rate function $V(r):=\frac{\|r\|^2}{2}$ (this is a standard application of the G\"artner Ellis Theorem);
	therefore, by a simple application of the contraction principle, we get the desired LDP holds with 
    good rate function $I_{Z^{(L+1)}(\underline{x})}$ defined by \eqref{eq:finitedimensional-rf} because the function
	$$(g^{(L)},r)=((g_{\alpha\beta}^{(L)})_{\alpha,\beta\in A},(r_{\alpha h})_{\alpha h})
	\mapsto g^{(L),\#}r=\big(\sum_{\gamma\in A}g_{\alpha\gamma}^{(L),\#}r_{\gamma h}\big)_{\alpha h}$$
	is continuous (here we also take into account Theorem X.1.1 in \cite{Bhatia}).
	
	We start with the case $L=1$. In this case the ratio $\frac{n_\ell(n)}{n_{\widehat{\ell}}(n)}$ in Condition \ref{cond:on-n1-nL}
	is trivially equal to 1 (because we have $\ell=\widehat{\ell}=1$); thus $v_n^*=n_1(n)$ and $\gamma_L=\gamma_1=1<\infty$.
	Firstly, it is easy to check that, by \eqref{eq:random-variances} with $\ell=1$, for all 
	$\eta\in\mathbb{R}^{|A|\times|A|}$ we have
	$$\lim_{n\to\infty}\frac{1}{v_n^*}\log\mathbb{E}[e^{v_n^*\langle\eta,G_n^{(1)}(\underline{x})\rangle}]
	=\langle\eta,C_b\mathbf{1}\rangle
	+\kappa(C_W\eta;g^{(0)}(\underline{x}))=:\Psi(\eta;g^{(0)}(\underline{x})).$$
	So, by the G\"artner Ellis Theorem on $\mathbb{R}^{|A|\times|A|}$ (note that $\Omega_2\subset\mathbb{R}^{|A|\times|A|}$), 
	we prove $(\bullet)$ for $L=1$ (so we necessarily have $\gamma_1=1$) noting that, by \eqref{eq:finitedimensional-rf-auxiliary-bis} with $\gamma_1=1$,
	$$\sup_{\eta\in\mathbb{R}^{|A|\times|A|}}\{\langle\eta,g^{(1)}\rangle-\Psi(\eta;g^{(0)}(\underline{x}))\}=J(g^{(1)}|g^{(0)}(\underline{x}))$$
	coincides with $I_{G^{(1)}(\underline{x})}(g^{(1)})$ (see \eqref{eq:finitedimensional-rf-auxiliary} for $L=1$, with a slight abuse of notation; indeed the
	infimum in \eqref{eq:finitedimensional-rf-auxiliary} disappears because $g^{(1)}$ is fixed, we have the unique constraint $g^{(0)}=g^{(0)}(\underline{x})$
	because the constraint $g^{(1)},\ldots,g^{(L-1)}\in\mathcal{S}_{|A|,C_b}$ is empty, and the sum in \eqref{eq:finitedimensional-rf-auxiliary} is reduced to single
	summand). We also remark that $I_{G^{(1)}(\underline{x})}$ is a good rate function; indeed (here we restrict to the case $\gamma_1=1$ for what we have said 
	above, but this restriction is not necessary) the function $\Psi(\cdot;g^{(0)}(\underline{x}))$ assumes finite values because $\sigma(\cdot)$ is bounded (by Condition \ref{cond:on-sigma}) and we can refer to Lemma 2.2.20 in \cite{DemboZeitouni}.
	
	Now we consider the inductive hypothesis, i.e. we assume that $(\bullet)$ holds for $L-1$ (for $L\geq 2$). In 
	what follows we prove that $(\bullet)$ holds by a suitable application of Proposition \ref{prop:Chaganty} with $\Omega_1=\Omega_2=\mathcal{S}_{|A|,C_b}$. So we have to check the three following items:
	\begin{enumerate}
		\item $I_{G^{(L)}(\underline{x})}$ is a good rate function;
		\item if we take $g_n^{(L-1)}\to g^{(L-1)}$ as $n\to\infty$ in $\Omega_1$, then the sequence of conditional distributions
		$$\{P(G_n^{(L)}(\underline{x})\in\cdot|G_n^{(L-1)}(\underline{x})=g_n^{(L-1)})\}_n$$
		satisfies the LDP on $\Omega_2$, with speed $v_n^*$ and rate function $J(\cdot|g^{(L-1)})$ defined by \eqref{eq:finitedimensional-rf-auxiliary-bis};
		\item the function $(g^{(L-1)},g^{(L)})\mapsto J(g^{(L)}|g^{(L-1)})$ is lower semicontinuous.
	\end{enumerate}
	Indeed, if these conditions hold, we have the following equality
	$$I_{G^{(L)}(\underline{x})}(g^{(L)}):=\inf_{g^{(L-1)}\in\mathcal{S}_{|A|,C_b}}\{J(g^{(L)}|g^{(L-1)})+I_{G^{(L-1)}(\underline{x})}(g^{(L-1)})\},$$
	which meets the expression in \eqref{eq:finitedimensional-rf-auxiliary}.
	
	We start with item 2. If we prove it for $\gamma_L<\infty$, then the proof for $\gamma_L=\infty$ is a consequence of the 
	application of Lemma \ref{lem:two-speed-functions} with
	$\{\psi_n\}_n=\{P(G_n^{(L)}(\underline{x})\in\cdot|G_n^{(L-1)}(\underline{x})=g_n^{(L-1)})\}_n$, 
	$H=J(\cdot|g^{(L-1)})$ in \eqref{eq:finitedimensional-rf-auxiliary-bis} which uniquely vanishes at 
	$$r_0=C_b\mathbf{1}+C_Wy(g^{(L-1)})\in\mathcal{S}_{|A|,C_b}$$
	(here we refer to $y(q)$ defined in \eqref{eq:unique-zero-LT}), $s_n=n_L(n)$ and $t_n=v_n^*$ (note that 
	$\frac{s_n}{t_n}=\frac{n_L(n)}{v_n^*}\to\gamma_L=\infty$). So, in what follows, we restrict our attention to the case 
	$\gamma_L<\infty$. We start noting that, by Lemma \ref{lem:equality-in-law} (and in particular
	eq. \eqref{eq:random-variances}), for all $\eta\in\mathbb{R}^{|A|\times|A|}$ we have
    $$\mathbb{E}[e^{\langle\eta,G_n^{(L)}(\underline{x})\rangle}|G_n^{(L-1)}(\underline{x})=g_n^{(L-1)}]
    =\exp\big(\langle\eta,C_b\mathbf{1}\rangle+n_L(n)\kappa(\frac{C_W}{n_L(n)}\eta;g_n^{(L-1)})\big);$$
	then, by Condition \ref{cond:on-sigma} (indeed $\kappa(\cdot;\cdot)$ assumes finite values because $\sigma(\cdot)$ is 
	bounded, and it is also a continuous function as discussed just after eq. \eqref{eq:def-kappa}), we get
	$$\displaylines{\lim_{n\to\infty}\frac{1}{v_n^*}\log\mathbb{E}[e^{v_n^*\langle\eta,G_n^{(L)}(\underline{x})\rangle}|G_n^{(L-1)}(\underline{x})
		=g_n^{(L-1)}]\phantom{oooooooooooooollllllllloooo}\cr\phantom{ooooooooooooooooooooooo}
	=\langle\eta,C_b\mathbf{1}\rangle+\gamma_L\kappa(\frac{C_W}{\gamma_L}\eta;g^{(L-1)})=:\Psi(\eta;g^{(L-1)}).}$$
	So, by the G\"artner Ellis Theorem on $\mathbb{R}^{|A|\times|A|}$ (note that $\Omega_2\subset\mathbb{R}^{|A|\times|A|}$), 
	we prove item 2 for $L$ with $\gamma_L<\infty$ noting that,
	by \eqref{eq:finitedimensional-rf-auxiliary-bis},
	$$\sup_{\eta\in\mathbb{R}^{|A|\times|A|}}\{\langle\eta,g^{(L)}\rangle-\Psi(\eta;g^{(L-1)})\}=J(g^{(L)}|g^{(L-1)}).$$
	In particular we can say that $J(\cdot|g^{(L-1)})$ is a good rate function, indeed $\Psi(\cdot;g^{(L-1)})$ assumes finite 
	values because $\sigma(\cdot)$ is bounded (by Condition \ref{cond:on-sigma}), and we can refer again to Lemma 2.2.20 in 
	\cite{DemboZeitouni} as we did above.
	
	For item 3 we have to check that, if $(g_n^{(L-1)},g_n^{(L)})\to (g^{(L-1)},g^{(L)})$ as $n\to\infty$ in $\Omega_1\times\Omega_2$, then
	\begin{equation}\label{eq:lsc-condition}
		\liminf_{n\to\infty}J(g_n^{(L)}|g_n^{(L-1)})\geq J(g^{(L)}|g^{(L-1)}).
	\end{equation}
	In order to do that (in both cases $\gamma_L<\infty$ and $\gamma_L=\infty$) one can check that, by \eqref{eq:finitedimensional-rf-auxiliary-bis}, 
	for all $\eta\in\mathbb{R}^{|A|\times|A|}$ we have
	$$J(g_n^{(L)}|g_n^{(L-1)})\geq\langle\eta,g_n^{(L)}\rangle-\Psi(\eta;g_n^{(L-1)}),$$
	where (here we take into account that, for the case $\gamma_L=\infty$, the function in \eqref{eq:finitedimensional-rf-auxiliary-bis} is the
	Legendre transform of a suitable linear function of $\eta$)
	$$\Psi(\eta;g^{(L-1)}):=\left\{\begin{array}{ll}
		\langle\eta,C_b\mathbf{1}\rangle+\gamma_L\kappa(\frac{C_W}{\gamma_L}\eta;g^{(L-1)})&\ \mbox{if}\ \gamma_L<\infty\\
		\langle\eta,C_b\mathbf{1}+C_Wy(g^{(L-1)})\rangle&\ \mbox{if}\ \gamma_L=\infty
	\end{array}\right.$$
	(actually the definition of $\Psi(\eta;g^{(L-1)})$ for $\gamma_L<\infty$ was already given above when we have checked item 2); then we can say 
	that \eqref{eq:lsc-condition} holds by letting $n$ go to infinity and by taking the supremum with respect to $\eta$.
	
	We conclude with item 1. In what follows we do not distinguish the cases $\gamma_L<\infty$ and $\gamma_L=\infty$. We refer to the final 
	statement of Proposition \ref{prop:Chaganty} with $\Omega_1=\Omega_2=\mathcal{S}_{|A|,C_b}$, and to the inductive hypothesis ($I_{G^{(L-1)}(\underline{x})}$ is a good rate function); then it suffices to show that
	$$(g^{(L-1)},g^{(L)})\mapsto J(g^{(L)}|g^{(L-1)})+I_{G^{(L-1)}(\underline{x})}(g^{(L-1)})$$
	is a good rate function (then we get the goodness of $I_{G^{(L)}(\underline{x})}$ by an application of the final statement of 
	Proposition \ref{prop:Chaganty}). Moreover, by Lemma \ref{lem:Chaganty}, it is enough to show that, for every compact $K\subset\Omega_1$ 
	and $a\geq 0$,
	$$\mathcal{U}_{K,a}:=\bigcup_{g^{(L-1)}\in K}\{g^{(L)}\in\Omega_2:J(g^{(L)}|g^{(L-1)})\leq a\}\ \mbox{is a compact set of}\ \Omega_2.$$
	We take a sequence $\{g_n\}_n$ in $\mathcal{U}_{K,a}$ and we have to show that there exists a subsequence that converges to a point 
	in $\mathcal{U}_{K,a}$. Firstly, for every $n$, there exists $h_n\in K$ such that $J(g_n|h_n)\leq a$. Then there exists a subsequence of 
	$\{h_n\}_n$ (which we call again $\{h_n\}_n$) that converges to some $\widehat{h}\in K$ (because $K$ is compact); moreover, for the 
	corresponding subsequence of $\{g_n\}_n$ (which we call again $\{g_n\}_n$), we have
	$$J(g_n|h_n)=\sup_{\eta\in\mathbb{R}^{|A|\times|A|}}\{\eta g_n-\Psi(\eta;h_n)\}.$$
	
	We remark that, since $\sigma(\cdot)$ is bounded by Condition \ref{cond:on-sigma}, we have $J(g|h)<\infty$ if
	$$C_b\leq g_{\alpha\beta}\leq C_b+C_W\|\sigma\|_\infty^2\ (\mbox{for all}\ \alpha,\beta\in A)$$
	for every $g=(g_{\alpha\beta})_{\alpha,\beta\in A}$ and $h$ (in particular see also \eqref{eq:random-variances}).
	Then, since $J(g_n|h_n)\leq a$, if we consider the notation $g_n=(g_{n;\alpha\beta})_{\alpha,\beta\in A}$, we have
	$$C_b\leq g_{n;\alpha\beta}\leq C_b+C_W\|\sigma\|_\infty^2\ (\mbox{for all}\ \alpha,\beta\in A).$$
	Thus there exists a compact set $\widetilde{K}\subset\mathcal{S}_{|A|,C_b}$ (which does not depend on $n$) such that 
	$g_n\in\widetilde{K}$; so there exists a subsequence of $\{g_n\}_n$ (which we call again $\{g_n\}_n$) which converges to 
	some $\widehat{g}\in\widetilde{K}$. In conclusion, by the item 3 checked above, we have
	$$a\geq\liminf_{n\to\infty}J(g_n|h_n)\geq J(\widehat{g}|\widehat{h});$$
	thus $\widehat{g}\in\mathcal{U}_{K,a}$ because $\widehat{g}\in\{g\in\Omega_2:J(g|\widehat{h})\leq a\}$ and $\widehat{h}\in K$.
	
	\begin{remark}\label{rem:more-advanced-result-for-random-covariances}
		Actually we can say that $\{(G_n^{(1)}(\underline{x}),\ldots,G_n^{(L)}(\underline{x}))\}_n$ satisfies the LDP on 
		$(\mathcal{S}_{|A|,C_b})^L$ with speed $v_n^*$ and good rate function $I_{G^{(1:L)}(\underline{x})}$ defined by
		$$I_{G^{(1:L)}(\underline{x})}(g^{(1)},\ldots,g^{(L)})=\sum_{\ell=1}^LJ(g^{(\ell)}|g^{(\ell-1)}),$$
		where $g^{(0)}=g^{(0)}(\underline{x})$ is defined by \eqref{eq:initial-variance}. This LDP yields the one stated in $(\bullet)$.
	\end{remark}

\subsection{Proof of Theorem \ref{th:finitedimensional-MD}}
	We follow the same lines of the proof of Theorem \ref{th:finitedimensional-LD} (we refer again to Lemma \ref{lem:equality-in-law}).
	For every $L\geq 1$ we take into account the statement $(\bullet)$ in the proof of Theorem \ref{th:finitedimensional-LD}
	together with Lemma \ref{lem:two-speed-functions} with $\{\psi_n\}_n=\{P(G_n^{(L)}(\underline{x})\in\cdot)\}_n$, $H=I_{G^{(L)}(\underline{x})}$
	in \eqref{eq:finitedimensional-rf-auxiliary} which uniquely vanishes at $r_0=\widehat{g}_{\underline{x}}^{(L)}$, $s_n=v_n^*$ and $t_n=1/a_n$
	(note that $\frac{s_n}{t_n}=a_nv_n^*\to\infty$). Then we have
	$$(\bullet\bullet):\ \left\{\begin{array}{l}
		\mbox{$\{G_n^{(L)}(\underline{x})\}_n$ satisfies the LDP on $\mathcal{S}_{|A|,C_b}$,}\\
		\mbox{with speed $1/a_n$ and good rate function $\widetilde{I}_{G^{(L)}(\underline{x})}=\Delta(\cdot;\widehat{g}_{\underline{x}}^{(L)})$.}
	\end{array}\right.$$
	So, by $(\bullet\bullet)$ and by the LDP of $\{(\sqrt{a_n}N_{\alpha h}^{(L+1)})_{\alpha h}\}_n$
	(which can be obtained by a standard application of the G\"artner Ellis Theorem), a simple application of 
	the contraction principle (already explained in the proof of Theorem \ref{th:finitedimensional-LD}) yields the desired LDP
    with rate function $\widetilde{I}_{Z^{(L+1)}(\underline{x})}$ defined by
	$$\widetilde{I}_{Z^{(L+1)}(\underline{x})}(z)=\inf_{g^{(L)}\in\mathcal{S}_{|A|,C_b},r\in\mathbb{R}^{|A|\times n_{L+1}}}
	\{\Delta(g^{(L)};\widehat{g}_{\underline{x}}^{(L)})+\|r\|^2/2:\widehat{g}^{(L),\#}r=z\}\quad 
	(z\in\mathbb{R}^{|A|\times n_{L+1}}),$$
	which coincides with the rate function in \eqref{eq:finitedimensional-rf-MD}.

\section{Large and moderate deviations of deep neural networks with ReLU pre-activation and single input}\label{sec:ReLUoneinput}
In this section we assume that $|A|=1$ (as in Remark \ref{rem:more-explicit-|A|=1})
and some notation can be simplified. In particular, for a standard Normal random variable $N$, we have 
$$\kappa(\eta;q):=\log\mathbb{E}[e^{\eta\sigma^2(\sqrt{q}N)}]\quad(\eta\in\mathbb{R}),$$
where $q\in\mathcal{S}_{|A|,C_b}=[C_b,\infty)$. Moreover we still consider the Legendre transform of $\kappa(\cdot;q)$ in \eqref{eq:def-kappa-star},
i.e.
$$\kappa^*(y;q)=\sup_{\eta\in\mathbb{R}}\{\eta y-\kappa(\eta;q)\}.$$

In this section, motivated by the literature, we consider the ReLU networks, i.e.
$$\sigma(x):=\max\{x,0\}.$$
In this case the function $\sigma(\cdot)$ is continuous and unbounded; therefore we cannot refer to Theorems 
\ref{th:finitedimensional-LD} and \ref{th:finitedimensional-MD}. However we can obtain the same results by considering suitable 
modifications of some parts of the proofs, which will be discussed in this section.

\subsection{Modifications of the proofs of Theorems \ref{th:finitedimensional-LD} and \ref{th:finitedimensional-MD}}\label{sub:4.1}
The parts with $\gamma_L=\infty$ still work well (and therefore the proof of Theorem \ref{th:finitedimensional-MD}).
On the other hand, when $\gamma_L<\infty$, we have to change some parts in which we use the G\"artner Ellis Theorem
in the proof of Theorem \ref{th:finitedimensional-LD}. For $q>0$ we have
\begin{multline*}
	\kappa(\eta;q)=\log\mathbb{E}[e^{\eta\sigma^2(\sqrt{q}N)}]=\log\mathbb{E}[e^{\eta qN^2\cdot 1_{\{N\geq 0\}}}]\\
	=\log\big(\mathbb{E}[e^{\eta qN^2}1_{\{N\geq 0\}}]+P(N<0)\big)
	=\log\big(\int_0^\infty e^{\eta qx^2}\frac{e^{-\frac{x^2}{2}}}{\sqrt{2\pi}}dx+1/2\big)\\
	=\log\big((\int_{-\infty}^\infty e^{\eta qx^2}\frac{e^{-\frac{x^2}{2}}}{\sqrt{2\pi}}dx+1)/2\big)
	=\left\{\begin{array}{ll}
		\log\frac{(1-2\eta q)^{-1/2}+1}{2}&\ \mbox{if}\ \eta<(2q)^{-1}\\
		\infty&\ \mbox{otherwise}.
	\end{array}\right.
\end{multline*}
Obviously this formula also holds for $q=0$ (this can happen when $C_b=0$), and we have $\kappa(\eta;0)=0$ for all 
$\eta\in\mathbb{R}$. Then, for every $q\geq C_b$, $\kappa(\eta;q)$ is finite in a neighborhood of the origin $\eta=0$, 
essentially smooth and lower semi-continuous (note that, under these conditions, the rate function $\kappa^*(\cdot;q)$ is 
good by Lemma 2.2.20 in \cite{DemboZeitouni}). Moreover we need to have
\begin{equation}\label{eq:pointwise-limit}
	\lim_{n\to\infty}\kappa(\eta_n;q_n)=\kappa(\eta;q),\ \mbox{as}\ q_n\to q\geq C_b\ \mbox{and}\ \eta_n\to\eta\in\mathbb{R};
\end{equation}
indeed one can easily check this condition (it is useful to distinguish the cases $\eta<(2q)^{-1}$, $\eta>(2q)^{-1}$ and 
$\eta=(2q)^{-1}$; if $q=0$ we always have $\eta<(2q)^{-1}$). We remark, when we presented the proofs of
Theorems \ref{th:finitedimensional-LD} and \ref{th:finitedimensional-MD}, the function $\sigma(\cdot)$ is bounded, 
and condition \eqref{eq:pointwise-limit} holds; indeed, in that case, we have a sequence of finite-valued functions, and the 
limit of this sequence is a finite-valued function.

Other modifications concern the details on how to check item 1 in the proof of Theorem \ref{th:finitedimensional-LD}.
We recall that we have (here we take into account that $|A|=1$)
$$J(g|h)=\sup_{\eta\in\mathbb{R}}\{\eta g-\Psi(\eta;h)\},$$
where
$$\Psi(\eta;h):=\left\{\begin{array}{ll}
	C_b\eta+\gamma_L\kappa(\frac{C_W}{\gamma_L}\eta;h)&\ \mbox{if}\ \gamma_L<\infty\\
	(C_b+C_Wy(h))\eta&\ \mbox{if}\ \gamma_L=\infty.
\end{array}\right.$$
We follow the same lines of that part of the proof of Theorem \ref{th:finitedimensional-LD}: for a compact subset $K$ of 
$[C_b,\infty)$, we have $h_n\in K$ for every $n$, and $J(g_n|h_n)\leq a$; then we have to check that there exists a subsequence 
of $\{g_n\}_n$ (which we call again $\{g_n\}_n$) that belongs to a compact subset $\widetilde{K}$ (say) of $[C_b,\infty)$. The 
case $\gamma_L=\infty$ is trivial because, for a standard Normal distributed random variable $N$, we have
$J(g_n|h_n)\leq a$ if and only if
$$g_n=C_b+C_Wy(h_n)=C_b+C_W\mathbb{E}[(\max\{\sqrt{h_n}N,0\})^2]=C_b+C_Wh_n\mathbb{E}[N^21_{N\geq 0}].$$
So, in what follows, we discuss the case $\gamma_L<\infty$. If $h=0$ (this can happen if $C_b=0$) we have $\kappa(\eta;0)=0$ 
for all $\eta\in\mathbb{R}$ and
$$J(g|0)=\sup_{\eta\in\mathbb{R}}\{\eta (g-C_b)\}=\Delta(g;C_b);$$
if $h>0$ we have $\kappa(\eta;q)=\kappa(\eta q;1)$ and
\begin{multline*}
	J(g|h)=\sup_{\eta\in\mathbb{R}}\big\{\eta g-\big(C_b\eta+\gamma_L\kappa(C_W\eta/\gamma_L;h)\big)\big\}
	=\sup_{\eta\in\mathbb{R}}\big\{\eta(g-C_b)-\gamma_L\kappa(C_W\eta h/\gamma_L;1)\big\}\\
	=\sup_{\eta\in\mathbb{R}}\big\{\eta((g-C_b)/h+C_b)-
	(C_b\eta+\gamma_L\kappa(C_W\eta/\gamma_L;1))\big\}
	=J((g-C_b)/h+C_b;1).
\end{multline*}
Then $J(g_n|h_n)\leq a$ yields $\frac{g_n-C_b}{h_n}+C_b\in\{y\in\mathbb{R}:J(y|1)\leq a\}$ if $h_n>0$ (and $\{y\in\mathbb{R}:J(y|1)\leq a\}$
is compact set that does not depend on $n$), and $g_n=C_b$ if $h_n=0$; thus (here we do not distinguish the cases $h_n>0$ and $h_n=0$) there 
exists $M\geq C_b$ such that
$$C_b\leq g_n\leq (M-C_b)h_n+C_b.$$
So $g_n$ belongs to a suitable compact subset $\widetilde{K}$ because $h_n$ belongs to a compact subset $K$.
Then we can conclude as for the case in which $\sigma$ is bounded and continuous.

\subsection{An explicit expression of $\kappa^*(\cdot;q)$}\label{sub:4.2}
Let $\kappa^\prime(\eta;q)$ be the derivative with respect to $\eta$, i.e.
$$\kappa^\prime(\eta;q)=\frac{q(1-2\eta q)^{-3/2}}{(1-2\eta q)^{-1/2}+1},$$
and let $\eta=\eta_{y,q}$ be the unique solution of $\kappa^\prime(\eta;q)=y$. Then we have
$$\kappa^*(y;q)=\left\{\begin{array}{ll}
	\eta_{y,q}y-\kappa(\eta_{y,q};q)&\ \mbox{if}\ y>0\\
	-\lim_{\eta\to-\infty}\kappa(\eta;q)=\log 2&\ \mbox{if}\ y=0\\
	\infty&\ \mbox{if}\ y<0.
\end{array}\right.$$
Moreover, if we consider the function $f(z)=\frac{z^3}{z+1}$ (for $z>0$), and we denote its inverse by
$f^\leftarrow$, then we have
$$(1-2\eta_{y,q}q)^{-1/2}=f^\leftarrow(y/q)$$
which yields
$$\eta_{y,q}=\frac{1-(f^\leftarrow(y/q))^{-2}}{2q}.$$
Thus, for $y>0$, we have
$$\kappa^*(y;q)=\frac{1-(f^\leftarrow(y/q))^{-2}}{2q}y-\log\frac{f^\leftarrow(y/q)+1}{2}.$$
Finally we remark that it is possible to give an explicit expression of $f^\leftarrow(y/q)$ in terms of the 
Cardano's formula for cubic equations. So we have an explicit expression of $\kappa^*(y;q)$, for which in 
general only a variational formula is available.

\section{Results for shallow neural networks and their sensitivities}\label{sec:L=1}
In this section we consider shallow Gaussian neural networks, i.e., the model \eqref{eq:Hanin-model} with $L=1$. We already
remarked that, in such a case, Condition \ref{cond:on-n1-nL} always holds. Throughout this section we set $n_1(n)=n$, so 
that $v_n^*=n$; then we have
$$Z_h^{(2)}(x)=b_h^{(2)}+\frac{1}{\sqrt{n}}\sum_{j=1}^n\sqrt{C_W}\widehat{W}_{hj}^{(2)}\sigma\big(b_j^{(1)}
+\sum_{r=1}^{n_0}W_{jr}^{(1)}x_r\big)\quad(h=1,\ldots,n_2),$$
where the random variables $\widehat{W}_{ij}^{(2)}$ are standard Normal distributed. We remark that we have a sum of $n$ i.i.d.
random variables (with respect to $j$), which depend on $h=1,\ldots,n_2$ and $x\in T\subset\mathbb{R}^{n_0}$; this particular
feature allows to establish some more results than the ones presented in the previous sections. In particular, motivated by 
the interest of the sensitivities with respect to the input $x\in T$ (see, e.g., \cite{FHMNP}), we can present 
results for some derivatives.

The aim of this section is to study large and moderate deviations of the $\mathbb{R}^{|A|\times n_2}$-valued sequence of random 
variables (in the derivatives below we have $x=(x_1,\ldots,x_{n_0})$)
\begin{equation}\label{eq:derivatives}
	\big\{\big(\frac 1{ \sqrt{n}}\frac{\partial^{s_h} Z_h^{(2)}(x)}
	{\partial x_1^{s_h}}\big|_{x=x_\alpha}
	\big)_{\alpha h}\big\}_n.
\end{equation}
Here we assume to have a finite set of inputs (i.e. $T=\{x_\alpha\}_{\alpha\in A}$ for a finite set $A$) and we take 
$s_h\in\{0,1\}$, $h=1,\ldots,n_2$; so, in particular, we define
$\frac{\partial^{0} Z_h^{(2)}(x)}{\partial x_1^{0}}:=Z_h^{(2)}(x)$.
Moreover, except for the case $s_1=\ldots=s_{n_2}=0$, we assume that $\sigma$ is, almost everywhere, differentiable; instead, if 
all the $s_h$ are equal to $0$, $\sigma$ may be not continuous.

It is easy to check that \eqref{eq:derivatives} reads
$$\big\{\big(\frac{b_h^{(2)}}{\sqrt n}1_{\{s_h=0\}}+\frac{1}{n}\sum_{j=1}^n
F_{hj}(x_\alpha)\big)_{\alpha h}\big\}_n,$$
where, for $x\in\mathbb{R}^{n_0}$,
\begin{equation}\label{eq:def-iid-rvs}
	F_{hj}(x):=\sqrt{C_W}\widehat{W}_{hj}^{(2)}
	\sigma^{(s_h)}\big(b_j^{(1)}+\sum_{r=1}^{n_0}W_{jr}^{(1)}x_{r}\big)
	(W_{j1}^{(1)})^{s_h}\quad (\mbox{for}\ j=1,\ldots,n)
\end{equation}
are $n$ i.i.d. $\mathbb{R}$-valued random variables.

In view of Propositions \ref{prop:derivative-LD} and \ref{prop:derivative-MD} below, it is useful to introduce the 
following functions:
$$\mathcal{I}_h(y_h):=\left\{\begin{array}{ll}
	\frac{y_h^2}{2C_b}&\ \mbox{if}\ C_b\neq 0\ \mbox{and}\ s_h=0\\
	\Delta(y_h;0)&\ \mbox{otherwise}
\end{array}\right.$$
(for $h=1,\ldots,n_2$ and $y_h\in\mathbb{R}$), 
\begin{eqnarray*}
	\Upsilon(\theta;\underline{x})\!\!\!&:=&\log\mathbb{E}\big[\exp\big(\sum_{\alpha\in A}\sum_{h=1}^{n_2}
	\theta_{\alpha h}F_{h1}(x_\alpha)\big)\big]\\
	&=&\log\mathbb{E}\big[\exp\big(\frac{C_W}{2}\sum_{h=1}^{n_2}\big(\sum_{\alpha\in A}\theta_{\alpha h}
	\sigma^{(s_h)}\big(b_1^{(1)}+\sum_{r=1}^{n_0}W_{1r}^{(1)}x_{\alpha,r}\big)
	(W_{11}^{(1)})^{s_h}\big)^2\big)\big]
\end{eqnarray*}
and
\begin{eqnarray*}
	\widetilde{\Upsilon}(\theta;\underline{x})\!\!\!&:=&\frac{1}{2}\sum_{\alpha,\beta\in A}\sum_{h_1,h_2=1}^{n_2}
	\left.\frac{\partial^2\Upsilon(\theta;\underline{x})}{\partial\theta_{\alpha h_1}\partial\theta_{\beta h_2}}\right|_{\theta=0}
	\theta_{\alpha h_1}\theta_{\beta h_2}\\
	&=&\frac{1}{2}\sum_{\alpha,\beta\in A}\sum_{h=1}^{n_2} 
	C_W\mathbb{E}\big[\sigma^{(s_{h})}\big(\sum_{r=1}^{n_0}W_{1r}^{(1)}x_{\alpha,r}\big)\sigma^{(s_{h})}
	\big(\sum_{r=1}^{n_0}W_{1r}^{(1)}x_{\beta,r}\big)(W_{11}^{(1)})^{2s_h}\big] \theta_{\alpha h}\theta_{\beta h}
\end{eqnarray*}
(for $\theta=(\theta_{\alpha h})_{\alpha\in A,h=1,\ldots,n_2}\in\mathbb{R}^{|A|\times n_2}$).

\begin{proposition}\label{prop:derivative-LD}
	Assume that the function $\Upsilon(\theta;\underline{x})$ is finite in a neighborhood of the origin 
	$\theta=0\in\mathbb{R}^{|A|\times n_2}$. Then the sequence 
	$\big\{\big(\frac{1}{\sqrt{n}}\frac{\partial^{s_h} Z_h^{(2)}(x)}
	{\partial x_1^{s_h}}|_{x=x_\alpha}\big)_{\alpha h}\big\}_n$ satisfies the LDP on $\mathbb{R}^{|A|\times n_2}$, with 
	speed $n$ and good rate function $I_{\partial Z^{(2)}(\underline{x})}$ defined by
	$$I_{\partial Z^{(2)}(\underline{x})}(z):=
	\inf\big\{\sum_{h=1}^{n_2}\mathcal{I}_h(y_h)+\Upsilon^*(f;\underline{x}):y_h+f_{\alpha h}=z_{\alpha h},\ \mbox{for}\ (\alpha,h)\in A\times\{1,\ldots,n_2\}\big\},$$
	where 
	$$\Upsilon^*(f;\underline{x}):=\sup_{\theta\in\mathbb{R}^{|A|\times n_2}}\big\{\sum_{\alpha\in A}\sum_{h=1}^{n_2}
	f_{\alpha h}\theta_{\alpha h}-\Upsilon(\theta;\underline{x})\big\}.$$
\end{proposition}
\begin{proof}
	The result can be proved by combining the LDP of $\{(\frac{1}{\sqrt{n}}b_h^{(2)})_h\}_n$ on $\mathbb{R}^{n_2}$ with speed 
	$n$ and good rate function $\sum_{h=1}^{n_2}\mathcal{I}_h(y_h)$ (this follows from a standard application of
	the G\"artner Ellis Theorem, and the independence of $b_1^{(2)},\ldots,b_{n_2}^{(2)}$), the LDP of
	$\{(\frac{1}{n}\sum_{j=1}^nF_{hj}(x_\alpha))_{(\alpha,h)}\}_n$ on 
	$\mathbb{R}^{|A|\times n_2}$ with speed $n$ and good rate function $\Upsilon^*(f;\underline{x})$ (by Cramér Theorem;
	see e.g. Theorem 2.2.30 and the subsequent Remark (a) in \cite{DemboZeitouni}), and a suitable application of the contraction
	principle (because the function $((y_h)_h,(f_{\alpha h})_{(\alpha,h)})\mapsto (y_h+f_{\alpha h})_{(\alpha,h)}$ is continuous).
\end{proof}

\begin{proposition}\label{prop:derivative-MD}
	Assume that the function $\Upsilon(\theta;\underline{x})$ is finite in a neighborhood of the origin 
	$\theta=0\in\mathbb{R}^{|A|\times n_2}$. Then, for every sequence of positive numbers $\{a_n\}_n$ such 
	that \eqref{eq:MD-conditions} holds with $v_n^*=n$, the sequence 
	$\big\{\big(\sqrt{a_n}\frac{\partial^{s_h} Z_h^{(2)}(x)}
	{\partial x_1^{s_h}}|_{x=x_\alpha}\big)_{\alpha h}\big\}_n$ satisfies the LDP 
	on $\mathbb{R}^{|A|\times n_2}$, with speed $1/a_n$ and good rate function $\widetilde{I}_{\partial Z^{(2)}(\underline{x})}$ defined by
	$$\widetilde{I}_{\partial Z^{(2)}(\underline{x})}(z):=
	\inf\big\{\sum_{h=1}^{n_2}\mathcal{I}_h(y_h)+\widetilde{\Upsilon}^*(f;\underline{x}):y_h+f_{\alpha h}=z_{\alpha h},\ \mbox{for}\ (\alpha,h)\in A\times\{1,\ldots,n_2\}\big\},$$
	where 
	$$\widetilde{\Upsilon}^*(f;\underline{x}):=\sup_{\theta\in\mathbb{R}^{|A|\times n_2}}\big\{\sum_{\alpha\in A}\sum_{h=1}^{n_2}
	f_{\alpha h}\theta_{\alpha h}-\widetilde{\Upsilon}(\theta;\underline{x})\big\}.$$
\end{proposition}
\begin{proof}
	It is similar to the proof of the previous proposition.
	The result can be proved by combining the LDP of $\{(\sqrt{a_n}b_h^{(2)})_h\}_n$ on $\mathbb{R}^{n_2}$ with speed 
	$1/a_n$ and good rate function $\sum_{h=1}^{n_2}\mathcal{I}_h(y_h)$, 
	the LDP of $\{(\frac{1}{\sqrt{n/a_n}}\sum_{j=1}^nF_{hj}(x_\alpha))_{(\alpha,h)}\}_n$ on 
	$\mathbb{R}^{|A|\times n_2}$ with speed $1/a_n$ and good rate function $\widetilde{\Upsilon}^*(f;\underline{x})$ (by Theorem 3.7.1 in \cite{DemboZeitouni}), and the same application of the contraction principle in the proof of Proposition \ref{prop:derivative-LD}.
\end{proof}


Thus it is important to find conditions on $\sigma$ under which the function $\Upsilon(\cdot;\underline{x})$ is 
finite in a neighborhood of the origin. For this purpose we present the following proposition.

\begin{proposition}\label{prop:for-the-finiteness-around-the-origin}
	Assume that there exists $M>0$ such that
	\begin{equation}\label{eq:*}
		\max_{h=1,\ldots,n_2}\big|\sigma^{(s_h)}\big(b+\sum_{r=1}^{n_0}w_rx_{\alpha,r}\big)
		w_1^{s_h}\big|\leq M\big(1+\Big|b+\sum_{r=1}^{n_0}w_rx_{\alpha,r}\big|\big)
	\end{equation}
	for  $s_h\in\{0,1\}$, $\alpha\in A$ and for every  $b,w_1,\ldots,w_{n_0}\in\mathbb{R}$. Then 
	$\Upsilon(\cdot;\underline{x})$ is finite in a neighborhood of the origin.
\end{proposition}

\begin{proof}
	It is easy to check (by taking into account \eqref{eq:*}) that, for some $C>0$, we have	
	$$
	\sum_{h=1}^{n_2}\big(\sum_{\alpha\in A}\theta_{\alpha h}
	\sigma^{(s_h)}\big(b+\sum_{r=1}^{n_0}w_rx_{\alpha,r}\big)
	w_1^{s_h}\big)^2 \leq  C\sum_{h=1}^{n_2}\sum_{\alpha\in A}\theta_{\alpha h}^2
	\big(1+b^{2}+\sum_{r=1}^{n_0}w_r^{2}\big);
	$$
	thus 
	\begin{multline*}
		\Upsilon(\theta;\underline{x})=\log\mathbb{E}\big[\exp\big(\frac{C_W}{2}\sum_{h=1}^{n_2}\big(\sum_{\alpha\in A}\theta_{\alpha h}
		\sigma^{(s_h)}\big(b_1^{(1)}+\sum_{r=1}^{n_0}W_{1r}^{(1)}x_{\alpha,r}\big)
		(W_{11}^{(1)})^{s_h}\big)^2\big)\big]\\
		\leq\log\mathbb{E}\big[\exp\big(C\sum_{h=1}^{n_2}\sum_{\alpha\in A}\theta_{\alpha h}^2
		\big(1+(b_1^{(1)})^2+\sum_{r=1}^{n_0}(W_{1r}^{(1)})^{2}\big)\big)\big],
	\end{multline*}
	and the final expression is finite in a neighborhood of the origin since the random variables $b_1^{(1)}$ and $W_{1r}^{(1)}$
	for $r=1,\ldots, n_0$ are Gaussian distributed and independent. 
\end{proof}

One could try to consider a stronger version of \eqref{eq:*} in which one refers to derivatives of order higher than the
first one. In such a case we would have products of more than two independent Gaussian random variables, and it would not
be possible to have a finite function $\Upsilon(\cdot;\underline{x})$ in a neighborhood of the origin. In our opinion it 
would be possible to overcome this problem by considering a simplified model (for instance the case $C_b=0$).

We conclude with some examples.

\begin{example}\label{ex:1}
	If we take 
	$s_1=\cdots=s_{n_2}=0$ then condition \eqref{eq:*} is the sublinearity condition on $\sigma$, i.e. there exists $C>0$ such that
	\begin{equation}\label{eq:16112023I}
		|\sigma(x)|\leq C(1+|x|).
	\end{equation}
	This condition holds, for instance,  for every bounded function $\sigma$ (but, if $\sigma$ is also continuous, we can refer to the results
	in this paper with $L\geq 1$), the function $\sigma(x)=\max\{x,0\}$ (already studied in Section 
	\ref{sec:ReLUoneinput}) concerning the ReLU networks, and the SWISH function $\sigma(x)=\frac{x}{1+e^{-x}}$.
\end{example}

\begin{example}\label{ex:2}
	We take
		$s_1=s_2=\cdots=s_{n_2}=1$ and we assume that the a.e. first derivative $\sigma^\prime$
		is bounded, with essential supremum $\|\sigma^\prime\|_\infty$. This condition holds for $\sigma(x)=\sin x$, the sigmoid function
		$\sigma(x)=\frac{1}{1+e^{-x}}$, the softplus function $\sigma(x)=\log(1+e^{-x})$, and we can still take the functions
		$\sigma(x)=\max\{x,0\}$ and $\sigma(x)=\frac{x}{1+e^{-x}}$ as in the previous example.
		
\end{example}

\begin{example}\label{ex:3}
    We take $n_2=2$, $s_1=0$ and $s_2=1$ and 
		we assume that $\sigma$ satisfies the sublinearity condition \eqref{eq:16112023I}, and that 
		$\sigma^\prime$ is bounded. We can realize that \eqref{eq:*} holds following the lines of the previous examples.
\end{example}

\paragraph{Acknowledgements.} The authors thank two anonymous referees for the careful reading of the manuscript and 
for some useful comments.

\paragraph{Funding.} All the authors acknowledge the support of Indam-GNAMPA. C.M. and B.P. acknowledge the
support of MUR Excellence Department Project awarded to the Department of Mathematics, University of Rome Tor Vergata
(CUP E83C23000330006) and the support of University of Rome Tor Vergata (project "Asymptotic Properties in Probability" 
(CUP E83C22001780005)).


\bibliography{biblio}
\bibliographystyle{amsplain}

\end{document}